\documentclass{amsart}

\newtheorem{thm}{Theorem}[section]

\newtheorem{dfn} [thm] {Definition}
 \newtheorem{twofactor}[thm]{}
 \newtheorem{ex}[thm]{Example}
 \newtheorem{remexamples}[thm]{Remarks and examples}
 \newtheorem{prop}[thm]{Proposition}
 \newtheorem{convention}[thm]{Convention}
\newtheorem{cor}[thm]{Corollary}

\newcommand{\D}{\mathcal{D}}
\newcommand{\C}{\mathcal{C}}
 \newcommand{\K}{\mathcal{K}}
\newcommand{\N}{\mathcal{N}}
 \newcommand{\M}{\mathcal{M}}
\newcommand{\E}{\mathcal{E}}

 \newcommand{\To}{\longrightarrow}

\input{xy}
\usepackage{amssymb}
\input{diagxy}
\xyoption{all}

\begin{document}

\title [$\lambda$-presentable morphisms]{$\lambda$-presentable
morphisms, injectivity \\ and (weak) factorization systems}

\author [M. H\'{e}bert] {Michel H\'{e}bert$^\dag$}

\thanks{$^\dag$Department of Mathematics, The American University in Cairo, Box 2511, Cairo,
Egypt. mhebert@aucegypt.edu}

\date {September 8, 2005}



\begin{abstract}

We show that every $\lambda_m$-injectivity class (i.e., the class
of all the objects injective with respect to some class of
$\lambda$-presentable morphisms) is a weakly reflective
subcategory determined by a functorial weak factorization system
cofibrantly generated by a class of $\lambda$-presentable
morphisms. This was known for small-injectivity classes, and
referred to as the ``small object argument". An analogous result
is obtained for orthogonality classes and factorization systems,
where $\lambda$-filtered colimits play the role of the transfinite
compositions in the injectivity case. $\lambda$-presentable
morphisms are also used to organize and clarify some related
results (and their proofs), in particular on the existence of
enough injectives (resp. pure-injectives).

\end{abstract}

\maketitle
\section*{Introduction}\label{S:0}

It is well-known that for every (small) set $\N$ of morphisms in a
locally presentable categories $\C$, (Cof($\N), \N^\Box$) is a
weak factorization system, where the class Cof($\N$) of
\emph{cofibrations} of $\N$ is the class ReChPo($\N$) of all
retracts of transfinite compositions of pushouts of the members of
$\N$ (complete definitions are recalled below). One may refer to
this as the small object argument, originating in homotopy theory
in the 60's (see \cite{Bk}). As a consequence, every
small-injectivity class is weakly reflective (and determined by
the weak factorization above). This was recently generalized in
various directions, for example in \cite{AHRT} and in \cite{Ch}.

Our main result (Theorem \ref{factorizations} below) generalizes
the small object argument as follows. A
\emph{$\lambda$-presentable morphism} of $\C$ is just a $\lambda
$-presentable object of the comma category $(A \downarrow \C)$ for
some object $A$. We show that for any class $\M$ of
$\lambda$-presentable morphisms, the class $\K = \M^\vartriangle$
of all objects which are injective with respect to all members of
$\M$ is weakly reflective and determined by a weak factorization
system (ReCh($\N), \N^\Box$), for some class $\N$ of
$\lambda$-presentable morphisms. Analogous results are shown for
orthogonality classes and their associated reflective
subcategories and factorization systems, where $\lambda$-filtered
colimits play the role of the transfinite compositions in the
injectivity case.

Variants of the construction of the transfinite composition used
in the theorem will have other applications, for example to give
simple proofs of the existence of ``sufficiently many" pure
subobjects, and of ``enough" absolutely pure objects in $\C$. Much
of Section \ref{S:2} of the paper is devoted to the promotion of
the $\lambda$-presentable morphisms as a natural concept to
organize and simplify various results. For example, we show that
FP-injective modules are precisely the modules which are injective
with respect to finitely presentable monomorphisms, and also that
in a locally presentable category $\C$ with the transferability
property, an object $A$ is injective (with respect to all monos)
if and only if there exists $\gamma$ such that $A$ is injective
with respect to all $\gamma$-presentable monos and all
$\gamma$-pure monos. We finally note that the existence of enough
(pure-) injectives in certain categories rests on the fact that
there exists $\gamma$ such that all (pure) monos in those
categories are transfinite composition of $\gamma$-presentable
(pure) monos.


\section{Main definitions and results}

For basic definitions and results on locally presentable and
accessible categories, we refer the reader to \cite{AR}. For
convenience, we recall the following.

If $\lambda$ is a regular infinite cardinal, an object $A$ of a
category $\C$ is $\lambda$-\textit{presentable} if the hom-functor
$\C(A,-)\colon \C \To \mathbf{Set}$  preserves $\lambda$-filtered
colimits. Then, $\C$ is $\lambda$-\textit{accessible} if it has
all $\lambda$-filtered colimits, as well as a (small) set $S$ of
$\lambda$-presentable objects such that every object in $\C$ is
the $\lambda $-filtered colimit of a diagram with all its vertices
in $S$. Finally, $\C$ is \textit{locally}
$\lambda$-\textit{presentable} if it is $\lambda $-accessible and
cocomplete (or, equivalently, complete). We write \textit{finitely
presentable} for $\omega$-presentable. We recall from \cite{H1}
the following definition.

\begin{dfn}
\textup{A morphism $f \colon A \To B$ in a category $\C$ is called
$\lambda$-\textit{presentable} if it is a $\lambda $-presentable
object of the comma category $(A \downarrow \C)$. Given a class
$\N$ of morphisms, $\lambda$-$\N$ will denote the class of all
$\lambda$-presentable morphisms in $\N$}.
\end{dfn}

Note that the comma category $(A \downarrow \C)$ is locally
$\lambda$-presentable for every object $A$ in a locally
$\lambda$-presentable category $\C$ (\cite{H3}). This is used in
particular to show the following, which we will need throughout
the paper:

\begin{thm} \label{ThPO} (\cite{H3}) Let $\C$ be a locally
$\lambda$-presentable category. Then $f\colon A\To B$ is $\lambda
$-presentable if and only if there exists a pushout diagram
\newbox\anglebox
\setbox\anglebox=\hbox{\xy \POS(0,75)\ar@{-} (0,0) \ar@{-}
(75,75)\endxy}
 \def\angle{\copy\anglebox}
 $$\bfig
 \scalefactor{0.8}
 \square[C`D`A`B;```f]
 \place(400,100)[\angle]
 \efig$$
 where $C$ and $D$ are $\lambda$-presentable.\\
\end{thm}

\begin{remexamples} \label{l-pres}
\end{remexamples}

\begin{enumerate} \item[(1)] In finitary varieties, the theorem
above amounts to say that a morphism $f \colon A \To B$ is
finitely presentable when $f$ provides a way to ``present $B$" by
adding less than $\lambda$ generators and relations to some
presentation of $A$ (see \cite{H2}). Using this, one can see for
example that a morphism $f \colon A \To B$ in the category
\textbf{CRng} of commutative rings with unit is finitely
presentable if and only if it is (isomorphic to) the canonical
homomorphism $A \To A[x_{1},...,x_{n}]/(p_{1},...,p_{m})$ for some
polynomials $p_{1},...,p_{m}$ over the variables
$x_{1},...,x_{n}$. In other words, if and only if $B$ is a
finitely presented $A$-algebra, with $f$ as its structure
morphism. This coincides with the ``finitely presented"
homomorphisms sometimes met in the Commutative Rings literature
(see \cite{P} for example).

\item[(2)] One can also use Theorem \ref{ThPO} to show that in the
category \textbf{Mod}-$R$ of all right $R$-modules, $R$ a ring,
the embedding $A \hookrightarrow B$ of a submodule is finitely
presentable if and only if the quotient $B/A$ is a finitely
presentable module. Note however that a finitely presentable
morphism may not have a finitely presentable kernel.

\item[(3)] As shown in \cite{H2}, Theorem \ref{ThPO} is true ``up
to a retraction in $(A \downarrow \C)$" in all $\lambda$-acessible
categories $\C$ with pushouts. Actually, the proof there can be
easily adapted to formulate a version for ``multipushouts" (in the
sense of Diers' multicolimits: see \cite{AR}, for example). In
particular, and more explicitely, this implies that in a locally
$\lambda$-multipresentable category $\C$, a morphism $f \colon A
\To B$ is $\lambda$-presentable if and only if it is the retract
in $(A \downarrow \C)$ of a component of the multipushout of some
morphism with $\lambda$-presentable domain and codomain.
Consequently, in the category \textbf{Fld} of fields, the finitely
presentable morphisms are just the finitely generated extensions.
Much of what follows will hold in this extended context.
\end{enumerate}

\begin{convention}
 For the rest of the paper, unless otherwise
specified, $\C$ will be a locally $\lambda$-presentable category.
\end{convention}

We will use $\lambda$-presentable morphisms to clarify and
generalize some constructions and results related to
small-injectivity and small-orthogonality classes. The following
type of construction has been used at least from the 1960's. We
will use it throughout the paper, so we describe it in detail for
easy reference.

\begin{twofactor} \label{twofactor} \textup{\textbf{Two
factorizations}}
\end{twofactor}

Given a class $\N$ of morphisms in $\C$, recall that a morphism is
a \textit{transfinite composition of morphisms in $\N$}, if it is
the canonical morphism $F(0) \To \textup{colim}\,F$ of a functor
$F \colon \delta \To \C$, where $\delta$ is an ordinal (seen as a
well-ordered category), $(F(\beta) \To F(\beta+1))\in\N$ for every
$\beta < \delta$, and $F(\beta) =
\mathop{\textup{colim}}\limits_{\gamma < \beta} F(\gamma)$ for
every limit ordinal $\beta < \delta$.

Given $f \colon A \To B$, and $\N$ a class of
$\lambda$-presentable morphisms which is \emph{stable under
pushouts} (i.e., the pushout of a member of $\N$ along any
morphism is in $\N$), we will construct two different
factorizations $f = f_\lambda f^*$ of $f$, with $f^*$ a
transfinite composition of morphisms in $\N$.

Put $A = A_0$ and $f = f_0$. Given any $\alpha < \lambda$, let
$\mathcal{G}_\alpha$ be a squeleton of the category of all
factorizations $(h,q) \colon  A_{\alpha} \to^h X \to^q B$ of
$f_{\alpha}$ with $h \in \N$, where a morphism from $(h,q)$ to
$(h',q')$ is a morphism $r \colon X \To X'$ in $\C$ such that the
following commutes:
$$\bfig
\morphism(0,0)<500,300>[A_{\alpha}`X;h]
\morphism(0,0)|b|/@{>}/<500,-300>[A_{\alpha}`X';h']
\morphism(500,300)|b|/@{>}/<0,-600>[X`X';r]
\morphism(500,300)/@{>}/<500,-300>[X`B;q]
 \morphism (500,-300)|b|/@{>}/<500,300>[X'`B;q']
\efig$$

Let $G_{\alpha}$ be the set of objects of $\mathcal{G}_\alpha$. We
consider the cone $S_{\alpha} = ( h_q~|~(h_q,q) \in G_{\alpha} )$,
and the cone $L_{\alpha}$ of a representative set of the
isomorphic classes of the underlying set of $S_{\alpha}$. We
emphasize that if $(h_q,q)$ and $(h_{q'},q')$ are distinct objects
in $G_{\alpha}$, $h_q = h_{q'}$ but $q \neq q'$, then $h_q$ and
$h_{q'}$ are distinct in $S_{\alpha}$ but not in $L_{\alpha}$. In
the first case (called the \emph{strict} case), we take the
colimit (= multiple pushout) of the source $S_{\alpha}$, and in
the second (\emph{loose}) case, the colimit of $L_{\alpha}$. Note
that $S_{\alpha} = L_{\alpha}$ if $f$ is the unique morphism from
$A$ to the terminal object \textbf{1} of $\C$. In both cases,
denote by $A_{\alpha^+}$ the colimit object, $g_\alpha = g_{\alpha
{\alpha}^+} \colon A_{\alpha} \To A_{\alpha^+}$ the canonical
morphism, and $f_{\alpha^+} \colon A_{\alpha^+} \To B$ the induced
morphism:
$$\bfig
\morphism (0,0)<500,300>[A_{\alpha}`X;h_q]
\morphism(0,0)|b|/@{>}/<500,-300>[A_{\alpha}`X';h_{q'}]
\morphism(0,0)/@{.>}/<1000,0>[A_{\alpha}`A_{\alpha^+};g_\alpha]
\morphism(500,300)|b|/@{.>}/<500,-300>[X`A_{\alpha^+};]
 \morphism (500,-300)|a|/@{.>}/<500,300>[X'`A_{\alpha^+};]
 \morphism (1000,0)|m|/@{.>}/<800,0>[A_{\alpha^+}`B;f_{\alpha^+}]
 \morphism (500,300)|a|/@{>}/<1300,-300>[X`B;q]
 \morphism (500,-300)|b|/@{>}/<1300,300>[X'`B;q']
\efig$$

For $\gamma \geq \alpha$, we define $g_{\alpha \gamma}$ by
composition (taking the colimit of the chain when needed). We
denote $g_{\alpha \lambda} \colon A_{\alpha} \To A_{\lambda}$ by
$k_{\alpha}$, and put  $f^* = k_0$.
$$\bfig
\morphism (-1200,0)/-{}/<0,0>[(1)`;]

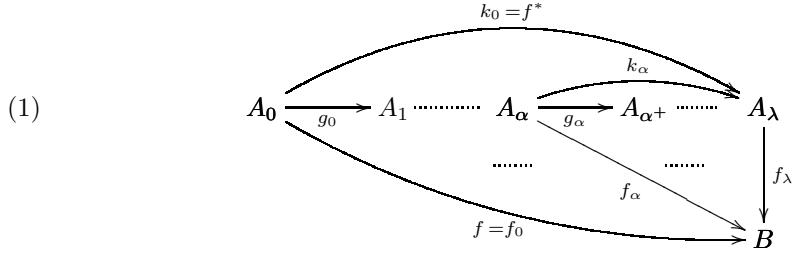
\morphism (-300,0)|b|<500,0>[A_0`A_1;g_0]

\morphism (-300,0)|a|/{@{>}@/^30pt/}/<1900,0>[A_0`A_\lambda
;k_0\,= f^*]

\morphism (-300,0)|b|/{@{>}@/^-15pt/}/<1900,-500>[A_0`B;f\,=f_0]

\morphism(250,0)/@{.}/<300,0>[`;]
\morphism(650,0)|b|<500,0>[A_{\alpha}`A_{\alpha^+};g_{\alpha}]

\morphism(650,0)|b|<950,-500>[A_{\alpha}`B;f_{\alpha}]

\morphism(550,-220)/@{.}/<200,0>[`;]
\morphism(1200,-220)/@{.}/<200,0>[`;]

\morphism(1150,0)/@{.}/<300,0>[A_{\alpha^+}`;]

\morphism(650,0)|a|/{@{>}@/^10pt/}/<950,0>[A_{\alpha}`A_{\lambda};k_{\alpha}]
\morphism(1600,0)|r|<0,-500>[A_{\lambda}`B;f_{\lambda}]

\morphism (2600,0)/-{}/<0,0>[`;]

 \efig$$
The factorization $f = f_{\lambda}f^*$ is the required one. Note
that the strict version of the construction is \emph{functorial}
in the following sense. Given $u \colon A \To A'$, $v \colon B \To
B'$ and $f'\colon A' \To B'$ such that $vf = f'u$, and
$f'_{\lambda}f'^*$ the strict factorization of $f'$, there is a
naturally induced $w \colon A_\lambda \To A'_\lambda$ such that
both squares in the following diagram commute:

$$\bfig

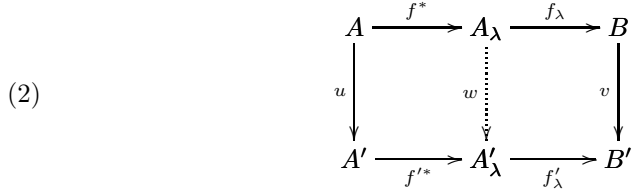
\morphism(-1250,-250)/-{}/<0,0>[(2)`;]
 \morphism (0,0)[A`A_\lambda;f^*]

\morphism (500,0)<500,0>[A_\lambda`B;f_\lambda]

\morphism (0,0)<0,-500>[A`A';u]

\morphism (1000,0)<0,-500>[B`B';v]

\morphism (0,-500)|b|[A'`A'_\lambda;f'^*]

\morphism (500,-500)|b|<500,0>[A'_\lambda`B';f'_\lambda]

 \morphism (500,0)/@{.>}/<0,-500>[A_\lambda`A'_\lambda;w]

\morphism (2500,-250)/-{}/<0,0>[`;]

 \efig$$

We spell out $g_0$ as a transfinite composition for future
reference, as we will consider interesting modifications of the
construction later. We well-order the underlying set $S^*_0 =
\{h_{\gamma}~|~\gamma < \delta = \delta(S_0) \}$ of $S_0$ (resp.
$L^*_0 = \{h_{\gamma}~|~\gamma < \delta = \delta(L_0) \}$ in the
loose case). Then we take the successive pushouts as follows.
First put $h_0 = h^1 = p^0$, and for $\gamma > 0$, let $p^\gamma$
be the pushout of $h_\gamma$ along $h^\gamma$, where $h^{\gamma
+1} = p^\gamma h^\gamma$, and for a limit ordinal $\gamma$,
$h^\gamma \colon A_0 \To X^\gamma$ is the canonical morphism to
the colimit $X^\gamma$ of the chain $(p^\beta ~|~ \beta <
\gamma)$. We illustrate for $\gamma \leq \omega$:
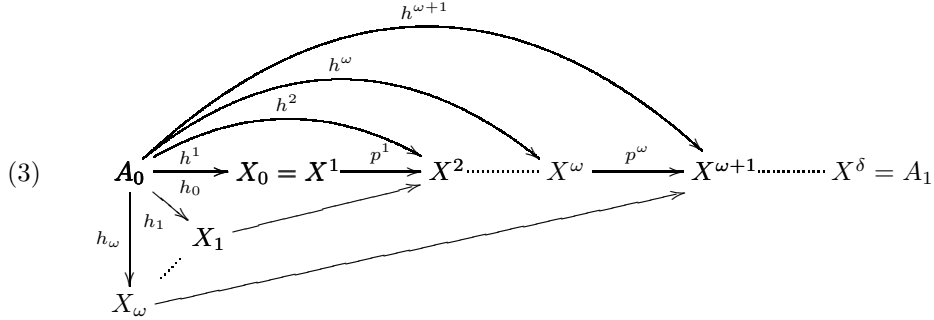
$$\bfig
\morphism(-400,0)/-{}/<0,0>[(3)`;]

 \morphism (0,0)|b|<600,0>[A_0`X_0 = X^1;h_0 ]

\morphism (0,0)|a|<600,0>[A_0`X_0 = X^1;h^1]

\morphism (0,0)|a|/{@{>}@/^20pt/}/<1200,0>[A_0`X^2;h^2]

\morphism (50,50)|a|/{@{>}@/^30pt/}/<1500,0>[` ;h^{\omega}]

\morphism
(0,0)|a|/{@{>}@/^55pt/}/<2250,0>[A_0`X^{\omega+1};h^{\omega+1}]

\morphism (800,0)<400,0>[`X^2;p^1]

\morphism (1250,0)/@{.}/<400,0>[`X^{\omega};]

\morphism (1750,0)<500,0>[`X^{\omega +1};p^\omega]

\morphism (120,-400)/@{.}/<100,100>[`;]

\morphism (2350,0)/@{.}/<500,0>[`X^{\delta}=A_1;]

\morphism (0,0)|b|<300,-250>[A_0`X_1;h_1]

 \morphism (0,0)<0,-500>[A_0`X_{\omega};h_{\omega}]

\morphism (100,-500)<2000,420>[`;]

 \morphism (300,-250)<800,200>[X_1`;]

 \morphism(3500,0)/-{}/<0,0>[`;]
 \efig$$

Clearly $p^{\gamma}$ is in $\N$ for all $\gamma < \delta$.\\

 We will first use \ref{twofactor} to obtain a generalization of the small
object argument. Before, we need to recall some definitions and
notations.

An object $K$ in $\C$ is \textit{injective with respect to}
(respectively \textit{orthogonal to}) a morphism $n \colon A \To
B$ if for all $u \colon A \To K$ there exists (resp. a unique) $w
\colon B \To K$ such that $wn = u$
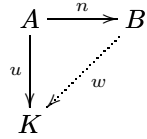
$$\bfig
 \scalefactor{0.8}
 \ptriangle/>`>`.>/[A`B`K;n`u`w]
 \efig$$

We write this as $n \triangledown K$ (resp. $n \perp K$). For a
class $\N$ of morphisms in $\C$, we define
$$\N^{\vartriangle} =\{K~| ~ n \triangledown K~ \textrm{for all}~ n\in\N  \}$$
$$\N^{\perp} =\{K~| ~ n \perp K ~\textrm{for all}~ n\in\N  \},$$
and for a class of objects $\K$, we define
$$\K^{\triangledown} =\{n~| ~ n \triangledown K ~\textrm{for all}~ K\in\K \}$$
$$\K^{\top} =\{n~| ~ n \perp K ~\textrm{for all}~ K\in\K \},$$
The class $\N^{\vartriangle}$ is also denoted by Inj$(\N)$ or
$\N$-Inj in the literature.

Given morphisms $p$ and $i$, we write $p$ $\Box$ $i$ (respectively
$p \perp i$) if for every commutative square $vp = iu$, there
exists a (resp. unique) morphism $w$ making both triangles commute
in the following diagram:
  $$\bfig
 \scalefactor{0.8}
 \square[ *`*`*`*;p`u`v`i]
 \morphism(500,500)|b|/-->/<-450,-450>[~`;w]
 \efig$$

If $\N$ is a class of morphisms, define the classes
$$^{\Box}\N =\{p~| ~ p ~\Box ~i ~ \textrm{ for all }   i \in
\N \}$$
$$\N^{\,\Box} =\{i~| ~ p ~\Box ~i ~ \textrm{ for all }   p \in
\N \}$$
$$\N^{\,\uparrow}=\{p~| ~ p \perp i ~ \textrm{ for all }   i \in
\N \}$$
$$\N^{\,\downarrow}=\{i~| ~p \perp i ~ \textrm{ for all }   p \in
\N \}.$$
 Then a \textit{weak factorization system} (resp.
a \textit{factorization system}) in $\C$ is a pair ($\E$,$\M$) of
classes of morphisms such that:
\begin{enumerate}
\item [(1)] $^{\Box}\M = \E$ (resp.
$\M^{\uparrow} = \E$),

\item [(2)] $\E^{\,\Box} = \M$ (resp. $\E^{\downarrow} = \M$), and

\item [(3)] every morphism $f$ in $\C$ has a
factorization $f = me$ with $m \in \M$ and $e\in \E$.
\end{enumerate}
(Note: These established notations are difficult to harmonize: the
tradition in the weak case has been to represent $p$ and $i$
vertically in the square above, so that what is oriented
``up/down" in factorization systems becomes ``left/right" in the weak ones.)\\

A class $\K$ of objects in $\C$ is a
\textit{$\lambda_{m}$-injectivity
  class} (resp. a \textit{$\lambda$-injectivity class}) if $\K =
  \N^{\vartriangle}$ for some class $\N$ of $\lambda$-presentable morphisms (resp.
  of morphisms with $\lambda$-presentable domains and codomains).
   \textit{$\lambda_{m}$-orthogonality
  classes} and \textit{$\lambda$-orthogonality classes} are defined
  similarly, replacing $\N^{\vartriangle}$ by $\N^\perp$.\\

We will use the following notations, given a class $\N$ of
morphisms:
\begin{enumerate}
\item [(1)] Ch$(\N$):= the class of all transfinite compositions
of morphisms in $\N$.

 \item [(2)] Po$(\N)$ := the
class of all pushouts of members of $\N$ (along any morphisms).

\item [(3)] $(A \downarrow \N)$ := the full subcategory of $(A
\downarrow \C)$, with $\N$ as its class of objects.

\item [(4)] Re($\N$) := the class of all retracts in $(A
\downarrow \C)$ of objects in $(A \downarrow \N)$.

\item [(5)] Cof$(\N)$ := ReChPo$(\N)$ (the class of the
\emph{cofibrations} of $\N$).

\item [(6)] Fc$_{\lambda}(\N)$ := the class of all the canonical
morphisms $A \To \textup{colim} \,U_\D$, with $\D$ a
$\lambda$-filtered subcategory of $(A \downarrow \N)$, and $U_\D
\colon \D \To \C$ its forgetful functor (defined by $U_\D (n
\colon A \To C) = C$).

\end{enumerate}

The following extends results in \cite{Bk}, \cite{HAR} and
\cite{Co}.

\begin{thm} \label{factorizations} Let $\M$ be a class of
$\lambda$-presentable morphisms in $\C$. Let $\N =
\lambda$-($\M^{\vartriangle\triangledown}$), and $\N_1 =
\lambda$-($\M^{\top\perp}$). Then
\begin{enumerate}

\item [(a)] $(\textup{ReCh}(\N),\N^{\,\Box})$ is a (functorial)
weak factorization system, and

\item [(b)] $(\textup{Fc}_{\lambda}(\N_1), \N_1^{\,\downarrow})$
is a factorization system.

\end{enumerate}

\end{thm}
\begin{proof}
(a)  $\N$ is easily seen to be stable under pushouts. Given $f
\colon A \To B$, we apply the strict factorization $f =
f_{\lambda}f^*$ in \ref{twofactor} to $f$. Hence $f^* \in
\textup{Ch}(\N)$, and we now show that $f_{\lambda} \in
\N^{\,\Box}$.

Consider a commutative square
$$\bfig
\scalefactor{0.8}
 \square[X`Y`A_{\lambda}`B;n`u`v`f_{\lambda}]
 \efig$$
with $n \in \N$. By Theorem \ref{ThPO}, there exists a pushout
square
\newbox\anglebox
\setbox\anglebox=\hbox{\xy \POS(0,75)\ar@{-} (0,0) \ar@{-}
(75,75)\endxy}
 \def\angle{\copy\anglebox}
 $$\bfig
 \scalefactor{0.8}
 \square[C`D`X`Y;z`s`t`n]
 \place(400,100)[\angle]
 \efig$$
with $C$ and $D$ $\lambda$-presentable. Then $su$ factorizes
through one of the $k_\alpha$'s (refer to the diagram (1) in
\ref{twofactor}), $su = k_{\alpha}l$, and we let $(l',z')$ be the
pushout of $(l,z)$, $z' \colon A_{\alpha} \To P$. We prove that
$z'\in \N$.

First, $z'$ is $\lambda$-presentable, since it is the pushout of
$z$. Then consider $h \colon A_\alpha \To K$, with $K \in
\M^{\,\vartriangle}$. Because $k_\alpha$ is a transfinite
composition of morphisms in $\N$, there exists $h' \colon
A_\lambda \To K$ such that $h'k_\alpha = h$. Then there exists $p
\colon Y \To K$ with $pn = h'u$ (since $n \in \N$). We have now
$$ptz = pns = h'us = h'k_{\alpha}l = hl,$$ so that there is a
(unique) $q \colon P \To K$ such that $qz' = h$ (and $ql' = pt$).
This shows that $z'\in \N$.
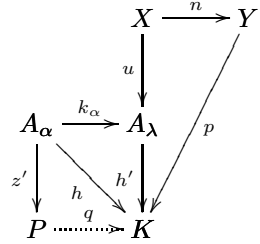
 $$\bfig
\scalefactor{0.8} \morphism (0,500)[X`Y;n]
 \morphism (0,500)<0,-500>[X`A_\lambda;u]
\morphism (500,500)|b|<-500,-1000>[Y`K;p] \morphism
(-500,0)[A_{\alpha}`A_\lambda;k_{\alpha}]
\morphism(-500,0)<0,-500>[A_\alpha`P;z']
\morphism(-500,0)|b|<500,-500>[A_\alpha`K;h]
\morphism(-500,-500)/@{.>}/<500,0>[P`K;q]
\morphism(0,0)<0,-500>[A_{\lambda}`K;h']
 \efig$$
Now, because $vtz = f_{\lambda}us = f_{\lambda}k_{\alpha}l$, the
pushout $l'z = z'l$ induces a (unique) morphism $r \colon P \To B$
such that $rz' = f_{\lambda}k_\alpha = f_\alpha$ (and $rl' = vt$).
Hence $A_\alpha \to^{z'} P \to^r B$ is in $G_\alpha$, and there
exists $z'' \colon P \To A_{{\alpha}^+}$ such that $z''z' =
g_\alpha$ and $f_{\alpha}z'' = r$.
\newbox\anglebox
\setbox\anglebox=
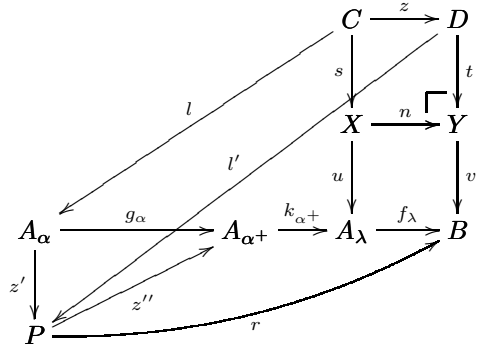
\hbox{\xy \POS(0,75)\ar@{-} (0,0) \ar@{-}
(75,75)\endxy}
 \def\angle{\copy\anglebox}
 $$\bfig
\scalefactor{0.8}
 \square |alra|[X`Y`A_{\lambda}`B;n`u`v`f_{\lambda}]
\square(0,500)[C`D`X`Y;z`s`t`]
 \place(400,600)[\angle]
 \morphism (-1500,0)<1000,0>[A_\alpha`A_{{\alpha}^+};g_{\alpha}]
\morphism (-500,0)[A_{{\alpha}^+}`A_\lambda;k_{{\alpha}^+}]
\morphism(-1500,0)<0,-500>[A_\alpha`P;z']
\morphism(0,1000)<-1500,-1000>[C`A_\alpha;l]
\morphism(500,1000)<-2000,-1500>[D`P;l']
\morphism(-1500,-500)|b|<1000,500>[P`A_{{\alpha}^+};z'']
\morphism(-1500,-500)|b|/{@/_12pt/}/<2000,500>[P`B;r]
 \efig$$
Then $us = k_{{\alpha}^+}g_{\alpha}l = k_{{\alpha}^+}z''z'l =
k_{{\alpha}^+}z''l'z$, and the pushout $ns = tz$ induces a
(unique) morphism $w \colon Y \To A_\lambda$ such that $wn = u$
and $wt = k_{{\alpha}^+}z''l'$. $w$ is the required diagonal: $v$
is the unique $x$ such that $xn = f_{\lambda}u$ and $xt = vt$, but
we have also $f_{\lambda}wn = f_{\lambda}u$ and $f_{\lambda}wt =
f_{\lambda}k_{{\alpha}^+}z''l' = f_{\alpha}z''l' = rl' = vt$. We
conclude that $f_{\lambda}w = v$, as needed.

The rest follows a known argument. For
$(\textup{ReCh}(\N),\N^{\,\Box})$ to be a weak factorization
system, what remains to be seen is that $\textup{ReCh}(\N) =\,
^\Box(\N^{\,\Box})$, since we will then have
$(\textup{ReCh}(\N))^\Box = (^\Box(\N^{\,\Box}))^\Box =
\N^{\,\Box}$). Its \emph{functoriality} refers to the property
with the same name mentioned in \ref{twofactor} (diagram (2)).

The inclusion $\textup{ReCh}(\N) \subseteq \,
 ^\Box(\N^{\,\Box})$ is clear, so let $(f \colon A \To B) \in \,
^\Box(\N^{\,\Box})$. We have seen above that $f$ factorizes as $A
\to^g C \to^h B$, with $h \in \N^{\,\Box}$ and $g \in
\textup{Ch}(\N)$. The commutative square $1_Bf = hg$ induces a
diagonal $d$ making the triangles commute in the following diagram
$$\bfig
\scalefactor{0.8}
 \square[A`B`C`B;f`g`1_B`h]
 \morphism (500,500) /{.>}/<-500,-500>[B`C;d]
 \efig$$
 so that $hd = 1_f$ in $(A \downarrow \C)$. Hence $hg = f$ is in
 $\textup{ReCh}(\N)$.

(b) Given $f \colon A \To B$, let  $\D$ be the subcategory of $(A
\downarrow \N_1)$ with its objects the morphisms appearing in the
cone $S_0$ (refer to \ref{twofactor}, with $\N = \N_1$), and its
morphisms those of the category $\mathcal{G}_0$. One can verify
that $\D$ is $\lambda$-filtered in $(A \downarrow \C$). Taking the
colimit of the forgetful functor $U_\D \colon \D \To \C$, the
canonical morphisms $A \to^h \textup{colim} \,
U_\D \to^g B$ gives
the required factorization in one step. The verification that the
induced morphism $g$ is in $\N^{\downarrow}$ is similar to the one
in the weak case, as well as the rest of the proof that
$(\textup{Fc}_{\lambda}(\N_1), \N_1^{\,\downarrow})$ is a
factorization system.
\end{proof}

\begin{remexamples} \label{Beke-Coste}
\end{remexamples}
\begin{enumerate}
\item [(1)] Let $\M$ be any (small) set of morphisms. We can
assume that members of $\M$ have $\lambda$-presentable domains and
codomains. Then Proposition 1.3 of \cite{Bk} says that
$(\textup{Cof}(\M),\M^{\,\Box})$ is a weak factorization system.
This amounts to replace
$\lambda$-($\M^{\vartriangle\triangledown}$) by Po$(\M)$, since
(Po($\M))^\Box$ = $\M^\Box$. Similarly, Corollary 3.3.4 of
\cite{Co} shows that $\lambda$-($\M^{\top\perp}$) can also be
replaced by Po$(\M)$ in this case (so that $(\textup{Fc}_\lambda
\textup{Po}(\M), \M^\downarrow)$ is a factorization system),
provided $\M$ admits a \emph{$\lambda$-strong left calculus of
fractions} (see II.2 in \cite{HAR}; note that the set of morphisms
in $\M^{\,\top\perp}$ which have $\lambda$-presentable domains and
codomains does admit a $\lambda$-strong left calculus of
fractions, as well as $\lambda$-($\N^{\top\perp}$) for any class
$\N$ of $\lambda$-presentable morphisms.)

Both facts are proved following the same line than in
\ref{factorizations}, but the proof is simpler in this case.

\item [(2)] Of course we have, in part (a) of the theorem, that
ReCh($\N$) = Cof($\N$) (since Po($\N$) = $\N$), so that the weak
factorization system is \emph{cofibrantly generated} by a class of
$\lambda$-presentable morphisms. An interesting problem would be
to find conditions under which it is cofibrantly generated by some
\emph{set}. More generally, one would like to be able to describe
ReCh($\lambda$-($\M^{\vartriangle\triangledown}$)) in a more
constructive way, from the elements of $\M$.\\
\end{enumerate}

A subcategory $\K$ of $\C$ is \textit{weakly reflective} if for
every $A\in\C$, there exists $r_A \colon A \To A^*$ in
$\K^{\triangledown}$ with $A^* \in \K$. If $\K$ is also closed
under retracts, we say it is \textit{almost reflective}. In
locally $\lambda$-presentable categories, we know that:
\begin{enumerate}
\item [(1)] (\cite{AR}) Reflective (resp. almost reflective)
subcategories are orthogonality (resp. injectivity) classes.
\item[(2)] (\cite{H2},\cite{H3}) $\lambda_{m}$-orthogonality
($\lambda_{m}$-injectivity) classes are reflective (almost
reflective). \end{enumerate}

In (2), the fact that a $\lambda_{m}$-injectivity class is almost
reflective follows from it being closed under products, because
its inclusion in $\C$ satisfies the Solution Set Condition (by
\cite{H3}, Lemma 4.2).

If $\M$ is a family of $\lambda$-presentable morphisms, then it is
easily seen that
$(\lambda$-$(\M^{\vartriangle\triangledown}))^\vartriangle =
\M^\vartriangle$  and $(\lambda$-$(\M^{\top\perp}))^\top =
\M^\top$. Since a morphism $B \To \textbf{1}$ to the terminal
object is in $\M^\downarrow$ iff $B \in \M^\top$ (resp. is in
$\M^\Box$ iff $B \in \M^\vartriangle$), the (weak) factorization
system of Theorem \ref{factorizations} determines a (almost)
reflective subcategory $\K = \M^\perp$ (resp.  $\K =
\M^\triangledown$), with (weak) reflectors $r_A \colon A \To R(A)$
from the appropriate factorization $A \to^{r_A} R(A) \to
\textbf{1}$. Concerning the weak case, as mentioned in
\ref{twofactor}, the strict and the loose factorizations of $A \to
\textbf{1}$ are the same, and $R$ and $r$ define respectively a
functor $R \colon \C \To \C$ and a natural transformation $r
\colon 1_{\C} \To R$. Note that, in contrast with the reflective
case, weak reflectors are generally not functorial in this sense
(see \cite{T} for more on this).

It is well-known that every reflective subcategory of a locally
presentable category is determined in this way by some
factorization system (see \cite{CHK}). Whether this holds for
almost reflective subcategories and weak factorization system
appears to be an open problem. However, we conclude from the
above:

\begin{cor} \label{Reflec-WeakFactor} Let $\alpha$ be a regular cardinal,
$\alpha \geq \lambda$. Then every almost reflective subcategory
$\K$ of $\C$ which is an $\alpha_m$-injectivity class is induced
by a functorial weak factorization system (Cof($\N), \N^{\,\Box})$
for some class $\N$ of $\alpha$-presentable morphisms (i.e., $\K =
\N^\vartriangle)$. (In particular, $\K$ is functorially almost
reflective).

\end{cor}

As for the reflective case, note that the proof of Theorem
\ref{factorizations}(b) also gives a construction of the
reflectors, extending \cite{HAR}, II.3.

\begin{prop} \label{reflections} Let $\M$ be a class of
$\lambda$-presentable morphisms, $\K = \M^{\vartriangle}$. Then
the following are equivalent:
\begin{enumerate}
\item [(i)] $\K$ is a $\lambda$-injectivity class,

\item [(ii)] $\K$ is closed under $\lambda$-filtered colimits,

\item [(iii)] every $m \in \M$ is the pushout of some morphism in
$\K^{\triangledown}$ with $\lambda$-presentable domains and
codomains.

\item [(iv)] $\K$ is determined by some (functorial) weak
factorization system $(\textup{Cof}(\N),\N^{\,\Box})$, for some
set $\N$ of morphisms with $\lambda$-presentable domains and
codomains.
\end{enumerate}

 The same is true for $\K = \M^{\perp}$, if one
replaces $\K^{\triangledown}$ by $\K^{\top}$, injectivity by
orthogonality, and the weak factorization system
$(\textup{Cof}(\N),\N^{\,\Box})$ by the factorization system
$(\textup{Fc}_\lambda \textup{Po}(\N), \N^\downarrow)$.
\end{prop}
\begin{proof}
(i) $\Rightarrow$ (ii) is straightforward. (ii) $\Rightarrow$
(iii) is Lemma 3.6 of \cite{H2}, together with Theorem \ref{ThPO}
above. (iii) $\Rightarrow$ (iv): let $\N$ be the set of all
morphisms in $\K^{\triangledown}$ with $\lambda$-presentable
domains and codomains, which have some pushout in $\M$. Then
$\N^{\vartriangle} = \M^{\vartriangle} = \K$. That
$(\textup{Cof}(\N),\N^{\,\Box})$ is a weak factorization system is
\cite{Bk}. (iv) $\Rightarrow$ (i) is trivial.

The orthogonality case is completely analogous (the argument in
the proof of Lemma 3.6 of \cite{H2} works the same).
\end{proof}

\section{Examples and applications}\label{S:2}

\begin{ex} \label{trivial} \textup{$\M$ = $\lambda$-Mor.}
\end{ex}

Let $\M$ be the class $\lambda$-Mor of all $\lambda$-presentable
morphisms in a locally $\lambda$-presentable category $\C$.
Obviously $\lambda$-$(\M^{\vartriangle\triangledown}) = \M$. Given
$f \colon A \To B$, we compare the two factorizations $f =
f_{\lambda}f^*$ in \ref{twofactor}.

In the strict case, we have $f_{\lambda} \in \M^{\,\Box}$ (Theorem
\ref{factorizations}), which mean that for every commutative
square
$$\bfig
 \scalefactor{0.8}
 \square[C`D`A_{\lambda}`B;z`s`t`f_{\lambda}]
 \efig$$
 where $z$ is $\lambda$-presentable, there exists $d \colon D
 \To A_{\lambda}$ such that $dz = s$ \emph{and} $f_{\lambda}d =
 t$ (This will actually force $f_{\lambda}$ to be an isomorphism: see
 below). Looking at the proof of \ref{factorizations}, one sees that
 in the case of the loose factorization, everything works the same, except that
the found diagonal $d$ will only be guaranteed to satisfy $dz =
s$. This means precisely that $f_{\lambda}$ is a
$\lambda$-\textit{pure mono} (by Theorem \ref{ThPO}, this is
equivalent to the definition in \cite{AR}: see also \cite{H1}).

The loose factorization does not lead to any weak factorization
system, but what is interesting in this case is the control that
one keeps on the presentability of the $A_\alpha$'s, because each
$L_{\alpha}$ is a cone made of (essentially) distinct
$\lambda$-presentable morphisms: using Proposition 2.3.11 of
\cite{MP} and Theorem \ref{ThPO}, we find easily a cardinal
$\gamma$, depending on $\C$ only, such that if $A$ is
$\beta$-presentable for some $\beta > (\gamma^{<\,\gamma})^+$,
then $A_{\gamma}$ is also $\beta$-presentable. Regarding $\C$ as a
locally $\gamma$-presentable category, and applying the loose
factorization to $f \colon A \To B$, then $A_{\gamma}$ is
$\beta$-presentable and $f_\gamma$ is a $\gamma$-pure mono, hence
$\lambda$-pure. This simplifies the proof of Theorem 2.33 of
\cite{AR} for the existence of ``sufficiently many" $\lambda$-pure
subobjects in $\C$. Note however that pushouts are needed in our
case, while Theorem 2.33 of \cite{AR} applies to all accessible
categories $\C$.

Back to the strict case, consider our $f_{\lambda} \in
\M^{\,\Box}$. As any morphism with domain $A_\lambda$, it is the
colimit of a $\lambda$-filtered diagram in $(A_\lambda \downarrow
\C)$ made of $\lambda$-presentable morphisms (see \cite{H3},
Proposition 2.6). From straightforward diagram chasing, we find
easily a right inverse to $f_{\lambda}$; since it is (pure) mono,
it is an isomorphism.

Hence the weak factorization system
$(\textup{ReCh}(\M),\M^{\,\Box})$ of the theorem is just the
trivial factorization system (Mor, Iso), and $\M^{\vartriangle}$
is $\C$ itself. Note that this implies in particular that every
morphism in $\C$ is a transfinite composition of
$\lambda$-presentable morphisms. We will see in \ref{injectives}
that the same property for (pure) monomorphisms is uncommon: it
will be used crucially to deduce the existence of ``enough (pure-)
injectives" in some categories.\\

\begin{ex} \label{Besserre} \textup{Integral closure.}
\end{ex}

In \cite{Be}, A. Besserre constructs an ``integral closure" of
rings, which can be described more easily by the construction
\ref{twofactor}. In $\C$ =\textbf{ CRng}, let $\M$ be the class of
all morphisms of the form $A \To A[x]/(p)$, with $p$ a monic
polynomial. Such a morphism is easily seen to be a finitely
presentable mono, actually the pushout of the canonical
homomorphism from the free ring on the set $\{a_0, ..., a_{n-1}\}$
of the coefficients of $p$, to the ring freely presented by the
set of generators $\{a_0, ..., a_{n-1}, x \}$ and the relation
$p(x) = 0$. By Proposition \ref{reflections}, $\M^{\vartriangle}$
is an $\omega$-injectivity class, and the construction in
\ref{twofactor} shows that the associated weak reflectors are
monos  (they actually have several other interesting properties:
see \cite{P}, where the members of $\M^{\vartriangle}$ are called
the \emph{absolutely integrally closed} rings).

\begin{ex} \label{Monos} \textup{$\M$ = $\lambda$-Mono.}
\end{ex}

Let $\M$ be the class $\lambda$-Mono of all $\lambda$-presentable
monomorphisms. For example, in the category \textbf{Mod}-$R$ of
all right $R$-modules, with $\lambda = \omega$,
$\M^{\vartriangle}$ is the class of the \emph{FP-injective}
modules: this can be seen by comparing Theorem 5.39 in \cite{NY}
with \ref{l-pres} (2) above. (The FP-injective modules are also
called \emph{absolutely pure} modules, but this terminology is
misleading in the present context: see below). More generally,
recall that an object in a locally $\lambda$-presentable category
is called $\lambda$\emph{-injective} (\cite{F}) if it is injective
with respect to all monomorphisms with $\lambda$-generated domain
and $\lambda$-presentable codomain (where $A$ is
\emph{$\lambda$-generated} if the hom-functor $\C(A,-)\colon \C
\To \mathbf{Set}$  preserves colimits of $\lambda$-filtered
diagrams of monomorphisms). We have:


\subsection*{Proposition} \emph{Let $\lambda_g$-Mono be the
class of all monomorphisms with $\lambda$-generated domain and
$\lambda$-presentable codomain in $\C$. Then
$$\lambda_g\textup{-Mono} \subseteq \lambda\textup{-Mono} \subseteq
\textup{Po}(\lambda_g\textup{-Mono}).$$ In particular, an object
$A$ is $\lambda$-injective iff it is injective with respect to all
$\lambda$-presentable monos.}

\begin{proof}
From \cite{AR}, the $\lambda$-generated objects are precisely the
strong quotients of the $\lambda$-presentable objects. Also, it is
easy to check that for any diagram $X \to^e Z \to^g Y$ with $e$ a
strong epi, $g$ is $\lambda$-presentable if $ge$ is. The first
inclusion then follows.

Let $f \colon A \To B$ be a $\lambda$-presentable monomorphism. By
Theorem \ref{ThPO}, there exists a pushout
\newbox\anglebox
\setbox\anglebox=\hbox{\xy \POS(0,75)\ar@{-} (0,0) \ar@{-}
(75,75)\endxy}
 \def\angle{\copy\anglebox}
 $$\bfig
 \scalefactor{0.8}
 \square[C`D`A`B;g`u`v`f]
 \place(400,100)[\angle]
 \efig$$
with $C$ and $D$ $\lambda$-presentable. Take the (StrongEpi, Mono)
factorization $C \to^e E \to^m D$ of $g$. Since $f$ is mono, there
exists a unique $d \colon E \To A$ making everything commute in
the obtained diagram. $E$ is $\lambda$-generated, and one verifies
easily that $f$ is the pushout of $m$ along $d$. The last
statement follows immediately.
\end{proof}

 Hence $(\lambda$-Mono)$^{\vartriangle}$ is the class of
all $\lambda$-injective objects. This means that it is actually a
small-injectivity class, since there is only (essentially) a set
of $\lambda$-generated objects in $\C$. (Actually, one cannot
expect to prove that a given $\lambda_m$-injectivity class is not
a small-injectivity class without using some large-cardinal
principle, as this would violates the Vopenka's principle: see
\cite{AR}, Theorem 6.27). From the proposition above, we have
$(\lambda_g\textup{-Mono})^\Box = (\lambda\textup{-Mono})^\Box$,
Po($\lambda_g$-Mono) = Po($\lambda$-Mono), and hence
$$(\textup{Cof}(\lambda_g\textup{-Mono}),(\lambda_g\textup{-Mono})^{\,\Box}) =
(\textup{Cof}(\lambda\textup{-Mono}),(\lambda\textup{-Mono})^{\,\Box}).$$

Note however that $(\lambda$-Mono)$^{\vartriangle}$ is not
necessarily a $\lambda$-injectivity class. For example, we will
see below that the class of all FP-injectives $R$-modules is an
$\omega$-injectivity class (in \textbf{Mod}-$R$) if and only if
$R$ is a coherent ring.

\subsection*{2.3.1}
Enough $\lambda$-injectives.

We say that $\C$ \emph{has enough $\lambda$-injectives} if each $A
\in \C$ is the subobject of some $\lambda$-injective. This is
easily seen to be equivalent to the weak reflectors $r_A \colon A
\To A_\lambda$ determined by the above weak factorization system
to be monos. We observe that if $\C$ has enough
$\lambda$-injectives, then
$\lambda\textup{-}((\lambda\textup{-Mono})^{\vartriangle\triangledown})
= \lambda\textup{-Mono}$, so that the weak factorization system
provided by Theorem \ref{factorizations} is just
(Cof($\lambda$-Mono),($\lambda$-Mono)$^{\,\Box})$ above: indeed,
$\lambda\textup{-}((\lambda\textup{-Mono})^{\vartriangle\triangledown})
\supseteq \lambda\textup{-Mono}$ is clear, and if $g \colon A \To
C$ is in
$\lambda\textup{-}((\lambda\textup{-Mono})^{\vartriangle\triangledown})$,
then it must be mono because $r_A \colon A \To A_\lambda$
factorizes through it.

Now, assume that every transfinite composition of monos in $\C$ is
mono (as in all locally \emph{finitely} presentable categories,
for example), and that $\C$ has the \emph{transferability
property}, i.e., the class Mono of all monomorphisms is stable
under pushouts (one easily shows that this is equivalent to
$\lambda$-Mono being stable under pushouts). Then, as shown in
\cite{F}, $\C$ has (functorially) enough $\lambda$-injectives: the
left-hand part $f^*$ of the
(Cof($\lambda$-Mono),($\lambda$-Mono)$^{\,\Box})$-factorization of
any morphism $f$ (= the strict factorization in \ref{twofactor}
with $\N = \lambda$-Mono) is clearly mono in this case. Actually,
since $\C$ is locally $\gamma$-presentable for any regular $\gamma
> \lambda$, one sees easily that for each $A \in \C$, we have a
transfinite chain $$A \To A_\lambda \To A_\gamma \To ...,$$ where
the composition $A \To A_\gamma$ is the weak reflector associated
with (Cof($\gamma$-Mono),($\gamma$-Mono)$^{\,\Box})$. Note that if
each such chain \emph{weakly stabilizes}, i.e., there exists
$\gamma$ such that $A_\gamma \To A_\zeta$ is a split mono for all
$\zeta
> \gamma$, then $\C$ has \emph{enough injectives}: this means that
each $A$ is the subobject of some \emph{injective} (= a member of
(Mono)$^\vartriangle$).

Of course, if $\C$ does not have these nice properties, we can
always apply Theorem \ref{factorizations} to $\N =
\lambda$-$(\M^{\vartriangle\triangledown})$, but the associated
weak reflectors might then not be monos.

Assume that $\C$ has enough $\lambda$-injectives. Now, if
$(\lambda$-Mono)$^{\vartriangle}$ is a $\lambda$-injectivity
class, Proposition \ref{reflections} implies in particular that
every monomorphism with $\lambda$-generated domain and
$\lambda$-presentable codomain is the pushout of some monomorphism
with $\lambda$-presentable domain and codomain. One can deduce
from that that every $\lambda$-generated subobject of a
$\lambda$-presentable object in $\C$ is $\lambda$-presentable.
This conclusion is also reached in Theorem 4-15 of \cite{F}, where
such categories are called locally $\lambda$-coherent. However the
category \textbf{Mod}-$R$ has the transferability property, and
hence has enough $\omega$-injectives, but it is locally finitely
coherent if only if $R$ is a coherent ring (see again \cite{F}).

\subsection*{2.3.2} Absolutely $\lambda$-pure objects.

An object $A$ is called \emph{absolutely $\lambda$-pure} if every
monomorphism from $A$ is $\lambda$-pure. A $\lambda$-injective
object is always absolutely $\lambda$-pure (use the fact that the
right-end morphism in the (StrongEpi, Mono)-factorization of a
$\lambda$-presentable morphism is also $\lambda$-presentable). We
note a simple but interesting consequence (in what follows, $A$ is
\emph{$\lambda$-pure-injective} if it is injective with respect to
all $\lambda$-pure monos):

\subsection*{Corollary} \emph{Suppose that $\C$ satisfies
the transferability property. Then an object in $\C$ is injective
if and only if it is $\lambda$-injective and
$\lambda$-pure-injective.}

\begin{proof}
 For the non-trivial direction, given $g \colon A \To M$ with
 $M$ $\lambda$-injective, and a mono $f
\colon A \To B$, the pushout $f'$ of $f$ along $g$ is mono, hence
$\lambda$-pure (since $M$ is absolutely $\lambda$-pure). If $M$ is
also $\lambda$-pure-injective, then $f'$ is a split mono, and the
result follows.
\end{proof}

An absolutely $\lambda$-pure object is not necessarily
$\lambda$-injective. More precisely, we know from \cite{F} that if
transfinite compositions of monos in $\C$ are mono, then the
following are equivalent: (i) $\C$ has enough
$\lambda$-injectives; (ii) the class $\lambda$-Mono is stable
under pushouts; (iii) $\C$ has the
 transferability property; (iv) the absolutely $\lambda$-pure objects are
$\lambda$-injective. Nevertheless, we have:

\subsection*{Proposition} \emph{Assume that transfinite
compositions of monos in $\C$ are mono. Then $\C$ has enough
absolutely $\lambda$-pure objects.}

\begin{proof}
The conclusion means that for each $A$, there exists a
monomorphism $A \To A^*$ with $A^*$ absolutely $\lambda$-pure. Our
proof is strongly inspired by the one of the existence of enough
existentially closed objects in \cite{F} (Theorem 6-3). However
the use of presentable morphisms will simplify it much. We apply a
modified version of \ref{twofactor} (with $\N = \lambda$-Mono and
$f \colon A \To \textbf{1}$), which seems to be of general
interest. To accommodate for the fact that $\N$ is not stable
under pushouts, we simply apply the rule: when taking the
successive pushouts from the well-ordered $L_\alpha$, just discard
the unpleasant results.

More precisely, and referring to the diagram (3) in
\ref{twofactor}, if the first pushout $p^1$ (of $h_1$ along $h^1$)
is not in $\N$ (i.e., is not a mono), replace it by the identity
on $X^1$ (and hence $h^2 = h^1$). Similarly for the pushout of
$h_2$ along $h^2$, etc. Then $f^* \colon A_0 \To A_\lambda$ is a
mono, and we show that $A_\lambda$ is absolutely $\lambda$-pure as
follows: given a commutative square
$$\bfig
\scalefactor{0.8}
 \square[C`D`A_{\lambda}`B;z`u`v`g]
 \efig$$
with $g$ mono, and $C$ and $D$ $\lambda$-presentable, we must find
a diagonal from $D$ to $A_\lambda$ making the upper triangle
commute (here we use the definition of purity in \cite{AR}). If $C
\to^m E \to^e D$ is the (StongEpi, Mono) factorization of $z$,
then there exists a unique $s \colon E \To A_\lambda$ making
everything commute. $E$ being $\lambda$-generated, $s$ must
factorize through one of the $k_\alpha$'s (see diagram (1) in
\ref{twofactor}), $s = k_{\alpha}l$, and we let $(l',z')$ be the
pushout of $(l,z)$, $z' \colon A_{\alpha} \To P$. Then $z'$ is
mono because $gk_\alpha$ factorizes through it. Also, the pushout
of $z'$ along any $h^\beta \colon A_\alpha \To X^\beta$ is a mono,
since $X^\beta \To B$ must factorize through it. This means that
$z'$ was not discarded, and it implies that $A_\alpha \To
A_{\alpha^+}$ must factorize through it. The morphism $D \to^{l'}
P \to A_\lambda$ is then the required diagonal.
\end{proof}

If $\C$ has the amalgamation property, i.e., the pushout of a mono
along a mono is a mono (the category of groups is an example), the
construction in the proof is just the (unmodified) one in
\ref{twofactor}. Note that the absolutely $\lambda$-pure objects
may still not form a weakly reflective subcategory of $\C$ in this
case, but $f^* \colon A \To A_\lambda$, where $f_\lambda f^*$ is
the factorization of $A \To \textbf{1}$, has the injectivity
property with respect to monos: for every \emph{mono} $g \colon A
\To C$ with $C$ absolutely $\lambda$-pure, there exists $h$ with
$hf^* = g$.

In \textbf{CRng} (which does not even have the amalgamation
property), a finitely presentable morphism $f \colon A \To B =
A[x_{1},...,x_{n}]/(p_{1},...,p_{m})$ is a mono if and only if the
set $\{p_{1},...,p_{m}\}$ of polynomials is ``consistent" over
$A$, in the sense that the $p_{i}$'s have a common root in
\emph{some} extension of $A$. Hence the construction above is a
kind of ``algebraic closure". But the term ``algebraically closed"
is confusing here: although a field can be seen to be
algebraically closed in the usual sense if and only if it is
absolutely pure in the category \textbf{CRng} (using the
Nullstellensatz, see \cite{F}), a field $K$ can be algebraically
closed in an extension field $L$ without the embedding being pure.
This simply means that there exists a polynomial on \emph{several}
variables with coefficients in $K$ which have a solution in $L$
but not in $K$ (see \cite{Po} for an example).


\begin{ex} \label{Pure} \textup{$\M_\delta =
\delta$-$_\omega$Pure.}
\end{ex}

For simplicity, we assume that $\C$ is locally \emph{finitely}
presentable. Given a regular infinite cardinal $\delta$, we let
let $\M_\delta = \delta$-$_\omega$Pure be the class of all
$\delta$-presentable monos which are ($\omega$-) pure (NOT to be
confused with the $\delta$-pure monos above). We show that there
exists $\gamma$ (depending on $\C$ only) such that $\C$ has
``enough" ($\delta$-$_\omega$Pure)-injectives for every $\delta
\geq \gamma$, i.e, such that for each $A$ there exists a pure mono
$A \To A^*$ with $A^*$ in $(\M_{\delta})^\vartriangle$.

The proof of Theorem 2.4 of \cite{BR} (referring back to
\cite{AR}), states that there exists $\gamma \geq \lambda$ such
that for any $\delta \geq \gamma$, each object $C$ of $\C$ is a
$\delta$-filtered colimit $(u_i \colon C_i \To C)$ with the
$C_i$'s $\delta$-presentable and the $u_i$'s pure monos. Regarding
$\C$ as locally $\delta$-presentable, it follows from the proof of
Theorem \ref{ThPO} (see Lemma 2.5 in \cite{H3}), that for every
$\delta$-presentable pure mono $f$, there exists a pushout diagram
\newbox\anglebox
\setbox\anglebox=\hbox{\xy \POS(0,75)\ar@{-} (0,0) \ar@{-}
(75,75)\endxy}
 \def\angle{\copy\anglebox}
 $$\bfig
 \scalefactor{0.8}
 \square[C`D`A`B;z`u`v`f]
 \place(400,100)[\angle]
 \efig$$
 where $C$ and $D$ are $\delta$-presentable and $u$ and $v$ are pure monos.
This implies that $vz$ is pure, and hence $z$ too. Since pure
monos are stable under pushouts (see \cite{AHT}), we have
$\delta$-$_\omega$Pure = Po($\N$) for the set $\N$ of the pure
monos with $\delta$-presentable domains and codomains. In
addition, the induced weak reflectors are pure monos, by their
construction (using \cite{H3}, 2.12(9)).

Note also that, just as in 2.3.1, it follows that the weak
factorization system provided by Theorem \ref{factorizations} is
just $$(\textup{ReCh}(\delta\textup{-}_\omega\textup{Pure}),
(\delta\textup{-}_\omega\textup{Pure})^{\,\Box}) =
(\textup{Cof}(\N),\N^{\,\Box}),$$ and that for each $A$ we have a
we have a transfinite chain $$A \To A_\delta \To A_\gamma \To
...,$$ ($\gamma > \delta$) where the composition $A \To A_\gamma$
is the (pure mono) weak reflector determined by the corresponding
weak factorization system. If each such chain weakly stabilizes,
then $\C$ has \emph{enough pure-injectives}: each $A$ is a pure
subobject of a pure-injective.


\begin{ex} \label{injectives} \textup{Categories with enough (pure-)
injectives.}
\end{ex}

Again we assume that $\C$ is locally finitely presentable. We
follow the line of the proof of Theorem 2.4 in \cite{BR}, to give
parallel proofs of the facts that the existence of \emph{effective
union} of subobjects (resp. pure subobjects) implies the existence
of enough (resp. pure-) injectives (\cite{BR} deals with the
pure-injectivity case). Nothing will be really new here, but the
use of $\gamma$-presentable morphisms allows to see interesting
connections, and will shorten the proof in \cite{BR}.

Suppose that there exists $\gamma$ such that $\C$ satisfies the
two conditions:
\begin{enumerate}
\item [(1)] $\C$ has enough $\gamma$-injectives (resp. enough
($\gamma$-$_\omega$Pure)-injectives), and

\item[(2)] every mono (resp. pure mono) in $\C$ is the transfinite
composition of $\gamma$-presentable monos (resp.
$\gamma$-presentable pure monos).

 \end{enumerate}

Then it is obvious from the definitions that $\C$ has enough
injectives (resp. pure-injectives) -- actually functorially so,
since ($\gamma$-Mono)$^\vartriangle$ and
($\gamma$-$_{\omega}$Pure)$^\vartriangle$ are
$\gamma_m$-injectivity classes.

From \ref{Pure}, $\C$ always satisfy (1) for the pure-injectivity
case, and from 2.3.1, it satisfies (1) for the injectivity case if
(and only if) is satisfies the transferability property. We don't
know if condition (2) is necessary for having enough injectives or
pure-injectives (or at least functorially so), but it seems to be
an interesting property to consider in its own right (see for
example \ref{trivial} and \ref{l-pres} above).

We now assume that subobjects (resp. pure subobjects) in $\C$ have
\emph{effective unions}, i.e., the induced morphism $P \to B$ from
the pushout $P$ of the pullback morphisms of a pair of (resp.
pure) monos with codomain $B$ is a (pure) mono. Typical examples
are the categories of modules.

From \cite{AR} (see also \ref{Pure} for the ``pure" case), we know
that in any locally finitely presentable category, there exists a
cardinal $\gamma$ such that:

(a) each object $C$ of $\C$ is a $\gamma$-filtered colimit $(u_i
\colon C_i \To C)$ with the $C_i$'s $\gamma$-presentable and the
$u_i$'s monos (resp. pure monos), and

(b) every subobject of a $\gamma$-presentable object has its domain
$\gamma$-presentable.

Let $S$ be the set of all (pure) subobjects of an object $B$,
partially ordered by transfinite compositions of
$\gamma$-presentable (pure) monos. Since $S$ is closed under
colimit of chains, any given (pure) mono $f \colon A \To B$ must
factorize $A \to^t C \to^g B$ through a maximal element $g \colon
C \To B$ of $S$, with $t$ a transfinite composition of (pure)
monos. We need to show that $g$ is iso.

If not, then there exists a (pure) mono $b_0 \colon B_0 \To B$,
with $B_0 ~\gamma$-presentable, which does not factorize through
$g$ (by (a)). We take the pullback $(c_0,g_0)$ of $(b_0,g)$. All
four morphisms are monos, and $g_0$ is $\gamma$-presentable
because its domain is (see (b) above).

Then we take the pushout $(c'_0,g'_0 \colon C \To P)$ of
$(c_0,g_0)$, and let $h \colon P \To B$ be the induced morphism.

$$\bfig
\scalefactor{1.3}
 \square[C_0`B_0`C`B;g_0`c_0`b_0`g]
\morphism(0,0)<300,200>[C`P;g'_0]
\morphism(500,500)<-200,-300>[B_0`P;c'_0]
\morphism(300,200)|b|/@{.>}/<200,-200>[P`B;h]
 \efig$$

Then $g'_0$ is a $\gamma$-presentable, being the pushout of a
$\gamma$-presentable morphism, and it is a (pure) mono, because
$g$ is. $h$ is a (pure) mono, by the existence of effective unions
of (pure) subobjects, hence $g'_0$ is an iso, by maximality of
$g$. This is a contradiction because $b_0$ now factorizes through
$g$.


\end{document}